\documentclass[journal]{IEEEtran}
\usepackage{amssymb}
\usepackage{epsfig}
\usepackage{graphicx}
\usepackage{colortbl}

\newtheorem{theorem}{Theorem}
\newtheorem{corollary}{Corollary}
\newtheorem{remark}{Remark}
\newtheorem{lemma}{Lemma}
\newtheorem{assumption}{Assumption}
\newtheorem{definition}{Definition}

\newcommand{\hs}{\hspace}

\ifCLASSINFOpdf
\else
\fi
\begin{document}
%
\title{Consensusability of Multi-agent Systems with Delay and
Packet Dropout Under Predictor-like Protocols}
%
%
%

\author{Juanjuan~Xu,~Huanshui~Zhang~and~Lihua Xie~
\thanks{*This work is supported by the Taishan Scholar Construction Engineering by Shandong Government, the National Natural
Science Foundation of China under Grants 61403235, 61573221, 61633014. }
\thanks{J. Xu is with School of Control Science and Engineering, Shandong University, Jinan, Shandong, P.R. China 250061.
        {\tt\small jnxujuanjuan@163.com}}
\thanks{H. Zhang is with School of Control
Science and Engineering, Shandong University, Jinan, Shandong, P.R.China 250061.
        {\tt\small hszhang@sdu.edu.cn}}
\thanks{L. Xie is with School of Electrical and Electronic Engineering, Nanyang Technological University, Nanyang Avenue,
Singapore 639798.
        {\tt\small ELHXIE@ntu.edu.sg}}%

}

\maketitle

\begin{abstract}
This paper considers the consensusability of multi-agent systems with delay and packet dropout.
By proposing a kind of predictor-like protocol, sufficient and necessary conditions are given for the mean-square
consensusability in terms of system matrices, time delay, communication graph and the packetdrop probability. Moreover, sufficient
and necessary conditions are also obtained for the formationability of multi-agent systems.

\end{abstract}


%
\IEEEpeerreviewmaketitle

\section{Introduction}

Multi-agent systems have attracted much attention in various scientific communities due to their
broad applications in many areas including distributed computation \cite{lynch}, formation control \cite{fax},
distributed sensor networks \cite{cortes}. Consensus is the most fundamental control problem in multi-agent
systems. Due to the fact that each individual agent lacks global knowledge of the whole
system and can only interact with its neighbors, one key issue of consensus is
to study conditions under which the consensus can be achieved under a given protocol and other is the design of a consensus protocol.
Numerous results have been reported in the literature for the design of distributed consensus protocols
for multi-agent systems. See \cite{shuailiu}, \cite{olfatisaber} and references therein.
For the consensusability problem, \cite{ma} and \cite{you} gave a necessary and sufficient condition for the continuous-time and discrete-time multi-agent systems
in the deterministic case respectively. \cite{zong} studied the case with multiplicative noise and time delay.

Time delays are unavoidable in information acquisition and transmission of practical multi-agent systems and should be taken
into account in designing the consensus protocol. An initial study is given in \cite{olfatisaber}
which provides a necessary and sufficient condition on the
upper bound of time delays under the assumption that all the delays
are equal and time-invariant. Sufficient conditions have been given in \cite{bliman} for average consensus with
constant, time varying and nonuniform time delays.
\cite{munz} sutdied the output consensus for multi-agent systems with
different types of time delays including communication delay, identical self-delay and different self-delay. \cite{cao} considered
discrete-time multi-agent systems with dynamically changing topologies and time-varying
communication delays.

On the other hand, random link failures or transmission
noises exist widely in networked multi-agent systems, which motivates the study of stochastic consensus problem.
In the literature, \cite{fagnani} provided two kinds of average consensus protocols which
are biased compensation method and balanced compensation method in the presence of random link failures.
It was shown in \cite{xiao} that the consensus value will diverge when the traditional consensus algorithms are applied in the presence of noises.
Under a fixed topology, necessary and sufficient conditions were given in \cite{li} for mean square average
consensus. \cite{huang} derived a sufficient condition for the switching
topologies case. For the multiplicative-noise case, \cite{li1} revealed that multiplicative noises may enhance the almost sure consensus, but may
have damaging effect on the mean square consensus.
\cite{liz} studied the mean square consensus for linear discrete-time systems by
solving a modified algebraic Riccati equation. \cite{zong1} considered the stochastic consensus
conditions. \cite{zong} gave the stochastic consentability analysis of linear multi-agent systems
with time delays and multiplicative noises. Though plenty works have been done for multi-agent systems with either time delay or multiplicative noise, there is little progress for discrete-time multi-agent systems with both input delay and packet dropout. The consensus problem for the latter remains challenging. Note that the optimal control problem for the single agent system case was only solved recently by \cite{zhang3}.

In this paper, we will study the consensusability problem of multi-agent systems with delay and packet dropout.
Different from the consensus protocols in the literature where the protocol is mostly in the feedback form of the current state or the delayed state and there exists a maximum delay within which consensus can be achieved,
a new kind of predictor-like consensus protocol is proposed in this paper to deal with the delay.
Sufficient and necessary conditions are given for the mean-square
consensusability in terms of system matrix, time delay, communication graph and the packet dropout probability under the predictor-like protocol. It will be shown that the derived results can be reduced to the deterministic case obtained in the literature.
Moreover, sufficient and necessary conditions are obtained for the formationability of multi-agent systems.

The remainder of the paper is organized as follows. Section II presents some preliminary knowledge about algebraic graph theory.
Problem formulation is given in Section III. Section IV shows preliminaries on modified Riccati equation. Main results are stated in Section V. Some concluding remarks are given in the last section. Related theorems and proofs are given in Appendix.

The following notation will be used throughout this paper: $R^n$
denotes the family of $n$-dimensional vectors; $x'$ denotes the
transpose of $x$; a symmetric matrix $M>0\ (\geq 0)$ means that $M$ is
strictly positive-definite (positive semi-definite). $\hat{x}(k|t)\triangleq
E[x(k)|\mathcal{F}_{t-1}]$ denotes the conditional expectation with
respect to the filtration $\mathcal{F}_{t-1}.$ $\lambda_i(A)$ means the $i$th eigenvalue of matrix $A.$

\section{Algebraic Graph Theory}
In this paper, the information exchange among agents is modeled by an undirected graph.
Let $\mathcal{G}=(\mathcal{V}, \mathcal{E}, \mathcal{A})$ be a diagraph with the set of vertices
$\mathcal{E}=\{1, \ldots, N\}$, the set of edges $\mathcal{E}\subset \mathcal{V}\times \mathcal{V}$,
and the weighted adjacency matrix $\mathcal{A}=[a_{ij}]\in \mathbf{R}^{N\times N}$ is symmetric.
In $\mathcal{G}$, the $i$-th vertex represents the $i$-th agent. Let $a_{ij}>0$ if and only if $(i, j)\in \mathcal{E}$,
i.e., there is a communication link between agents $i$ and $j$. Undirected graph $\mathcal{G}$ is connected
if any two distinct agents of $\mathcal{G}$ can be connected via a path that follows the edges of $\mathcal{G}$.
For agent $i$, the degree is defined as $d_{i}\triangleq \sum_{j=1}^{N}a_{ij}$. Diagonal matrix
$\mathcal{D}\triangleq diag\{d_{1}, \ldots, d_{N}\}$ is used to denote the degree matrix of diagraph $\mathcal{G}$. Denote
the Laplacian matrix by $L_G=\mathcal{D}-\mathcal{A}$. The eigenvalues of $L_G$ are denoted by $\lambda_i(L_G)\in R, i=1,\cdots, N,$ and an ascending order
in magnitude is written as $0=\lambda_1(L_G)\leq \cdots\leq\lambda_N(L_G),$ that is,  the Laplacian matrix $L_G$ of an undirected graph has at
least one zero eigenvalue and all the nonzero eigenvalues are in the
open right half plane. Furthermore, $L_G$ has exactly one zero eigenvalue
if and only if $G$ is connected\cite{godsil}.


\section{Problem Formulation}

Consider a multi-agent system as depicted in Fig. 1 where the dynamic is given by
\begin{eqnarray}
x_i(k+1)&=&Ax_i(k)+\gamma(k)Bu_i(k-d),\nonumber\\
&&~~~~~~~~~~~~~i=1,\cdots,N,\label{m8}
\end{eqnarray}
while $x_i\in R^n$ is the state of the $i$th agent, $u_i\in R^m$ is the control input of the $i$th agent, $A,B$ are
constant matrices with appropriate dimensions. $d$ represents the input delay. $\gamma(k) = 1$
denotes that the data packet has been successfully delivered to
the plant, and $\gamma(k)=0$ signifies the dropout of the data packet. Without loss of generality, the random
process $\{\gamma(k),k\geq0\}$ is modeled as an independent and identically
distributed (i.i.d.) Bernoulli process with probability distribution $P(\gamma(k) = 0) = p$ and $P(\gamma(k) = 1) = 1- p,$ where $p\in(0,1)$
is said to be the packet dropout rate. The initial values are given by $x_i(0), u_i(-1),\cdots, u_i(-d).$ Note that the channel fading and
time delay occur simultaneously due to the unreliable network placed in the path from the controller $i$ to the agent $i$. Moreover, the information
exchange between the controllers of agent $i$ and $j$ happens in the controller processor.

\begin{remark}
Noting that the random process $\gamma$ is identical for which we will derive some necessary and sufficient conditions for consensusability of multi-agent systems
with both delay and packet dropout. The derived results will provide insights into the interplay among system dynamic, delay and network topology and demonstrate the advantage of the predictor-like consensus protocol. They could also shed some light on resolving the non-identical $\gamma$ case which is interesting and is left for our future study.
\end{remark}

We further make the following general assumption.
\begin{assumption}\label{a1}
All the eigenvalues of $A$ are either on or outside the unit circle, $B$ has full column rank.
\end{assumption}

\begin{assumption}\label{a3}
System $(A, B,0,A^dB)$ is mean-square stabilizable, that is, for the system
\begin{eqnarray}
x(k+1)=Ax(k)+Bu(k)+\nu(k)A^dBu(k)\nonumber
\end{eqnarray}
where $\nu(k)$ is a sequence of white noise with zero mean and unit covariance, there exists a
feedback controller $u(k)=Kx(k)$ where $K$ is a time-invariant matrix such that
the closed-loop system is mean-square stable, i.e. $\lim_{k\rightarrow\infty}E\|x(k)\|^2=0.$
\end{assumption}

\begin{assumption}\label{a2}
The undirected graph is connected.
\end{assumption}

\begin{center}
\begin{figure}\scalebox{0.35}{
\hs{0mm}\includegraphics{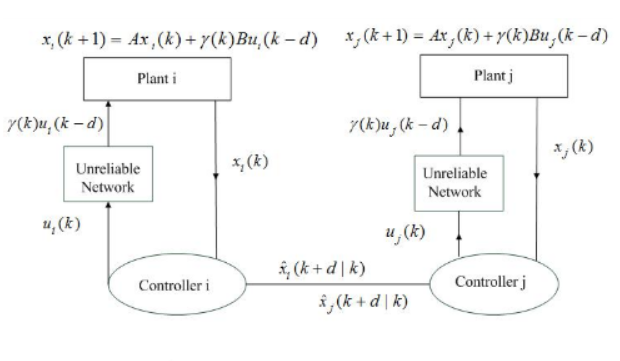}}
\caption{Multi-agent system with unreliable networks}
\end{figure}
\end{center}

Denote $w(k)=\gamma(k)-E[\gamma(k)],$ then system (\ref{m8}) is reformulated as
\begin{eqnarray}
x_i(k+1)&=&Ax_i(k)+(1-p)Bu_i(k-d)\nonumber\\
&&+w(k)Bu_i(k-d),~~~i=1,\cdots,N,\label{m18}
\end{eqnarray}
where $\{w(k), k\in N\}$ is a sequence of random variables defined
on $(\Omega, \mathcal{F}, \mathcal{P}; \mathcal{F}_k)$ with $E[w(k)] = 0$ and $E[w(k)w(s)] = p(1-p)\delta_{ks}.$
We simply denote $\mu=1-p$ and $\sigma^2=p(1-p).$

In the literature \cite{fax}, \cite{olfatisaber}, the relative state $x_j(k)-x_i(k)$ between agents is used to design the
consensus protocol like $u_i(k)=K\sum_{j\in N_i}\big[x_j(k)-x_i(k)\big].$ Differently in this paper, we firstly calculate the following predictor using
each agent's own state and historical inputs for $k\geq d$ in the way that
\begin{eqnarray}
&&\hat{x}_i(k|k-d)\nonumber\\
&=&E[x_i(k)|\mathcal{F}_{k-d-1}]\nonumber\\
&=&A^dx_{i}(k-d)+\mu\sum_{j=1}^{d}A^{j-1}Bu_i(k-d-j).
\end{eqnarray}
Then the relative predictor $\hat{x}_j(k|k-d)-\hat{x}_i(k|k-d)$ is applied to design the consensus protocol. To be specific,
the distributed protocol for $k\geq d$ is described as
\begin{eqnarray}
u_i(k-d)=K\sum_{j\in N_i}\big[\hat{x}_j(k|k-d)-\hat{x}_i(k|k-d)\big].\label{m19}
\end{eqnarray}
The aim is to find sufficient and necessary conditions for the mean-square consensusability of multi-agent system (\ref{m18}) under protocol
(\ref{m19}) where the definition on the mean square consensusability is given below.
\begin{definition}
The discrete-time multi-agent system (\ref{m18})
with a fixed undirected graph is said to be mean-square consensusable under protocol
(\ref{m19}) if for any finite initial values $x_i(0), u_i(-d),\cdots, u_i(-1),$ there exists a control gain $K$
such that the controller (\ref{m19}) enforces consensus, i.e.
$
\lim_{k\rightarrow\infty}E\|x_j(k)-x_i(k)\|^2=0,~\forall i, j=1,\cdots,N.$
\end{definition}

By substituting (\ref{m19}) into (\ref{m18}), the closed-loop multi-agent system becomes
\begin{eqnarray}
x_i(k+1)&=&Ax_i(k)+\mu BK\sum_{j\in N_i}\big[\hat{x}_j(k|k-d)\nonumber\\
&&-\hat{x}_i(k|k-d)\big]+BK\sum_{j\in N_i}\big[\hat{x}_j(k|k-d)\nonumber\\
&&-\hat{x}_i(k|k-d)\big]w(k),~~k\geq d.\label{m9}
\end{eqnarray}
Let $X(k)=\left[
           \begin{array}{ccc}
             x_1(k) & \cdots & x_N(k) \\
           \end{array}
         \right]',~\hat{X}(k|k-d)=\left[
           \begin{array}{ccc}
             \hat{x}_1(k|k-d) & \cdots & \hat{x}_N(k|k-d) \\
           \end{array}
         \right]',
$ then (\ref{m9}) can be reformulated as
\begin{eqnarray}
X(k+1)&=&(I_N\otimes A)X(k)-\mu(L_G\otimes BK)\hat{X}(k|k-d)\nonumber\\
&&- w(k)(L_G\otimes BK)\hat{X}(k|k-d),~~k\geq d.\label{m10}
\end{eqnarray}
Denote $\bar{X}(k)=\frac{1}{N}\sum_{i=1}^Nx_i(k),$ then
\begin{eqnarray}
\bar{X}(k+1)&=&\frac{1}{N}(1_N\otimes I_n)'X(k+1)\nonumber\\
&=&A\bar{X}(k)-(\mu/N)(1_N'L_G\otimes BK)\hat{X}(k|k-d)\nonumber\\
&&-(w(k)/N)(1_N'L_G\otimes BK)\hat{X}(k|k-d)\nonumber\\
&=&A\bar{X}(k),\label{m1}
\end{eqnarray}
where $1_N'L_G=0$ has been used in the derivation of the last equality.
Given the initial condition $\bar{X}(0)=\frac{1}{N}\sum_{i=1}^Nx_i(0)$ and equation (\ref{m1}),
it yields that $\bar{X}(k)$ is deterministic. This further implies that $\hat{\bar{X}}(k|s)=E[\bar{X}(k)|\mathcal{F}_{s-1}]=\bar{X}(k)$
for any positive integer $s.$ We now present the dynamic equation of $\delta(k+1)=X(k+1)-(1_N\otimes I_n)\bar{X}(k+1)$ with
$\hat{\delta}(k|k-d)=\hat{X}(k|k-d)-(1_N\otimes I_n)\hat{\bar{X}}(k|k-d)=\hat{X}(k|k-d)-(1_N\otimes I_n)\bar{X}(k).$
It is obtained by subtracting (\ref{m1}) from (\ref{m10}) that
\begin{eqnarray}
\delta(k+1)
&=&(I_N\otimes A)X(k)-\mu(L_G\otimes BK)\hat{X}(k|k-d)\nonumber\\
&&- w_k(L_G\otimes BK)\hat{X}(k|k-d)-(I_N\otimes A)\nonumber\\
&&\times (1_N\otimes I_n)\bar{X}(k)\nonumber\\
&=&(I_N\otimes A)X(k)-\mu(L_G\otimes BK)\hat{X}(k|k-d)\nonumber\\
&&- w_k(L_G\otimes BK)\hat{X}(k|k-d)-(I_N\otimes A)\nonumber\\
&&\times (1_N\otimes I_n)\bar{X}(k)+\mu(L_G\otimes BK)(1_N\otimes I_n)\nonumber\\
&&\times \hat{\bar{X}}(k|k-d)+w_k(L_G\otimes BK)(1_N\otimes I_n)\nonumber\\
&&\times \hat{\bar{X}}(k|k-d)\nonumber\\
&=&(I_N\otimes A)\delta(k)-\mu(L_G\otimes BK)\hat{\delta}(k|k-d)\nonumber\\
&&- w_k(L_G\otimes BK)\hat{\delta}(k|k-d),~~k\geq d.\label{m12}
\end{eqnarray}
Select $\phi_i\in R^N$ such that $\phi_i'L_G=\lambda_i(L_G)\phi_i'$ and form an unitary
matrix $\Phi=\left[
               \begin{array}{cccc}
                 \frac{1}{\sqrt{N}} & \phi_2 & \cdots & \phi_N \\
               \end{array}
             \right]
$ to transform $L_G$ into a diagonal form
\begin{eqnarray}
diag\{0,\lambda_2(L_G),\cdots,\lambda_N(L_G)\}=\Phi'L_G\Phi.\nonumber
\end{eqnarray}
Let $\tilde{\delta}(k)=(\Phi'\otimes I_n)\delta(k)=\left[
                                                                                       \begin{array}{ccc}
                                                                                         \tilde{\delta}_1(k) & \cdots & \tilde{\delta}_N(k) \\
                                                                                       \end{array}
                                                                                     \right].$
Together with the property of Kronecker product, it holds that $\tilde{\delta}_1(k)=0$
and for $k\geq d$ and $i=2,\ldots,N,$
\begin{eqnarray}
\tilde{\delta}_i(k+1)&=&A\tilde{\delta}_i(k)-\mu\lambda_i(L_G)BK\hat{\tilde{\delta}}_i(k|k-d)\nonumber\\
&&- \lambda_i(L_G)w_k BK\hat{\tilde{\delta}}_i(k|k-d).\label{m4}
\end{eqnarray}

\begin{theorem}\label{the3}
The multi-agent system (\ref{m18}) achieves mean-square consensus if and only if the systems in (\ref{m4}) are mean-square stable simultaneously.
\end{theorem}
\emph{Proof.} ``Necessity" The simultaneous mean-square stability of (\ref{m4}) follows from the derivation of (\ref{m9})-(\ref{m4}).

``Sufficiency" Since $\lim_{k\rightarrow\infty}E\|\tilde{\delta}_i(k)\|^2=0,$ then $\lim_{k\rightarrow\infty}E\|\delta_i(k)\|^2=0.$
This implies that $\lim_{k\rightarrow\infty}E\|x_i(k)-\bar{X}(k)\|^2=0$ for $i=1,\cdots,N.$ Thus,
\begin{eqnarray}
&&\lim_{k\rightarrow\infty}E\|x_j(k)-x_i(k)\|^2\nonumber\\
&\leq&\lim_{k\rightarrow\infty}E\|x_j(k)-\bar{X}(k)\|^2+
\lim_{k\rightarrow\infty}E\|x_i(k)-\bar{X}(k)\|^2\nonumber\\
&=&0.\nonumber
\end{eqnarray}
This gives the consensus of (\ref{m18}). The proof is now completed. \hfill $\blacksquare$

\section{Preliminaries on Modified Riccati Equation}

Based on Theorem \ref{the3}, the simultaneous stabilizability of the systems in (\ref{m4}) is necessary for consensusability.
To this end, we shall present some results with respect to the stabilizability criterion and further investigate a corresponding
modified algebraic Riccati equation. Firstly, the following equivalent conditions have been given in \cite{zhang1}.
\begin{lemma}\label{the1}
The following statements are equivalent.
\begin{enumerate}
  \item System
  \begin{eqnarray}
x(k+1)&=&Ax(k)+\mu Bu(k-d)\nonumber\\
&&+w_k Bu(k-d)\label{m7}
\end{eqnarray} is mean-square stable under the controller $u(k-d)=K\hat{x}(k|k-d).$
  \item System
  \begin{eqnarray}
x(k+1)=Ax(k)+\mu Bu(k)+w_k A^dBu(k)\label{m3}
\end{eqnarray}
is mean-square stabilizable under the controller $u(k)=Kx(k).$
  \item For any $Q>0,$ there exist matrices $K$ and $P >0$
satisfying the following equation:
\begin{eqnarray}
P&=&Q+(A+\mu BK)'P(A+\mu BK)\nonumber\\
&&+\sigma^2K'B'(A')^dPA^dBK.
\end{eqnarray}
\item There exist matrices $K$ and $P>0$
satisfying the following equation:
\begin{eqnarray}
P&>&(A+\mu BK)'P(A+\mu BK)\nonumber\\
&&+\sigma^2K'B'(A')^dPA^dBK.
\end{eqnarray}
\end{enumerate}
\end{lemma}
In particular, it has also been shown in \cite{zhang1} that the existence of a unique positive definite solution to
the algebraic Riccati equation
\begin{eqnarray}
P&=&A'PA+Q-\mu^2A'PB\Big[R+\mu^2B'PB\nonumber\\
&&+\sigma^2B'(A')^dPA^dB\Big]^{-1}B'PA\nonumber
\end{eqnarray}
is necessary and sufficient for the mean-square stabilizability of system (\ref{m7}) with $Q>0.$
Motivated by the results in \cite{zhang1}, we define the parameterized algebraic Riccati equation (PARE)
\begin{eqnarray}
P&=&A'PA+Q-\gamma A'PB\Big[R+B'PB\nonumber\\
&&+B'(A')^dPA^dB\Big]^{-1}B'PA\label{m22}
\end{eqnarray}
and denote
\begin{eqnarray}
g_{\gamma}(P)&=&A'PA+Q-\gamma A'PB\Big[R+B'PB\nonumber\\
&&+B'(A')^dPA^dB\Big]^{-1}B'PA,\label{r1}\\
\Phi(K,P)&=&(1-\gamma)(A'PA+Q)+\gamma(F_1'PF_1\nonumber\\
&&+F_2'PF_2+K'RK+Q),\label{r2}\\
\Psi(K,P)&=&F_1'PF_1+F_2'PF_2+K'RK+Q,\label{r3}
\end{eqnarray}
where $F_1=A+BK,F_2=A^dBK.$

\begin{theorem}\label{the10}
Consider the PARE (\ref{m22}). Let $A$ be unstable, $(A,B,0,A^dB)$ is mean-square stabilizable and $Q>0,R>0.$ Then the following hold.
\begin{enumerate}
  \item The PARE has a unique strictly positive definite solution if and only if $\gamma>\gamma_c,$ where $\gamma_c$ is the critical value defined as
\begin{eqnarray}
\gamma_c=\inf\{\gamma\in [0,1]|P=g_{\gamma}(P),P>0\}.\nonumber
\end{eqnarray}
  \item The critical value $\gamma_c$ satisfies the following analytical bounds:
\begin{eqnarray}
\underline{\gamma}\leq \gamma_c\leq \overline{\gamma}\nonumber
\end{eqnarray}
where $\underline{\gamma}$ and $\overline{\gamma}$ are defined by
\begin{eqnarray}
\underline{\gamma}&=&arginf_{\gamma}\{\exists S|(1-\gamma)A'SA+Q=S, S\geq0\}\nonumber\\
\overline{\gamma}&=&arginf_{\gamma}\{\exists (K,P)|P>\Phi(K,P)\}\nonumber
\end{eqnarray}
  \item The critical value can be numerically computed by the solution of the following quasiconvex LMI optimization problem
\begin{eqnarray}
&&\gamma_c=argmin _{\gamma}\Delta_\gamma(Y,Z)>0,0\leq Y\leq I\nonumber\\
&&\Delta_\gamma(Y,Z)=\nonumber\\
&&\left[
    \begin{array}{cccccc}
      Y & Y & \sqrt{\gamma}ZR^{\frac{1}{2}} & &  &  \\
      Y & Q^{-1} & 0 &  &  &  \\
      \sqrt{\gamma}R^{\frac{1}{2}}Z' & 0 & I &  &  &  \\
      \sqrt{\gamma}(AY+BZ') & 0 & Y  &  &  &  \\
      \sqrt{\gamma}A^dBZ' & 0 & 0 &  &  &  \\
      \sqrt{1-\gamma}AY & 0 & 0 &  &  &  \\
    \end{array}
  \right.\nonumber\\
&&\hspace{-12mm}\left.
    \begin{array}{cccccc}
       &  &  & \sqrt{\gamma}(AY+BZ')' & \sqrt{\gamma}(A^dBZ')' & \sqrt{1-\gamma}YA' \\
       &  &  & 0 & 0 & 0 \\
       &  &  & 0 & 0 & 0 \\
       &  &  & Y & 0 & 0 \\
       &  &  & 0 & Y & 0 \\
       &  &  & 0 & 0 & Y \\
    \end{array}
  \right]\nonumber
\end{eqnarray}
\end{enumerate}
\end{theorem}
\emph{Proof.} Based on Theorem \ref{the6}, \ref{the7} \ref{the8} and \ref{the9} in Appendix, the results follow by using similar proof to
that of Lemma 5.4 in \cite{sinopoli}.

\section{Mean-square Consensusability}

Denote for $i=2, \cdots,N,$
\begin{eqnarray}
\gamma_i=\frac{\mu^2}{\mu^2+\sigma^2}\frac{4\Big(\lambda_i(L_G)[\lambda_2(L_G)+\lambda_N(L_G)]-\lambda_i^2(L_G)\Big)}{[\lambda_N(L_G)+\lambda_2(L_G)]^2}.\nonumber
\end{eqnarray}

It is noted that
\begin{eqnarray}
\gamma_2
&=&\frac{\mu^2}{\mu^2+\sigma^2}\frac{4\lambda_2(L_G)\lambda_N(L_G)}{[\lambda_N(L_G)+\lambda_2(L_G)]^2}\nonumber\\
&=&\frac{\mu^2}{\mu^2+\sigma^2}\Big[1-\Big(\frac{\lambda_N(L_G)-\lambda_2(L_G)}{\lambda_N(L_G)+\lambda_2(L_G)}\Big)^2\Big].\nonumber
\end{eqnarray}

We now present the main result of the mean-square consensusability for multi-agent system (\ref{m18}).
\begin{theorem}\label{the2}
Let Assumption \ref{a1}-\ref{a2} hold.
If $\gamma_2>\gamma_c$ where $\gamma_c$ is given in Theorem \ref{the10},
then the multi-agent system (\ref{m18}) is mean-square consensusable under protocol
(\ref{m19}).
\end{theorem}
\emph{Proof.} Consider the Riccati equation
\begin{eqnarray}
P&=&A'PA+Q-\gamma_iA'PB\Big[R+B'PB\nonumber\\
&&+B'(A')^dPA^dB\Big]^{-1}B'PA.\label{m5}
\end{eqnarray}
Since \begin{eqnarray}
&&\frac{4\Big(\lambda_i(L_G)[\lambda_2(L_G)+\lambda_N(L_G)]-\lambda_i^2(L_G)\Big)}{[\lambda_N(L_G)+\lambda_2(L_G)]^2}\nonumber\\
&&-\frac{4\lambda_2(L_G)\lambda_N(L_G)}{[\lambda_N(L_G)+\lambda_2(L_G)]^2}\nonumber\\
&=&\frac{4\Big(\lambda_2(L_G)-\lambda_i(L_G)\Big)\Big(\lambda_i(L_G)-\lambda_N(L_G)\Big)}{[\lambda_N(L_G)+\lambda_2(L_G)]^2}\geq0,\nonumber
\end{eqnarray}
then it follows that $\gamma_i\geq \gamma_2>\gamma_c$ for $i>2.$ Using Theorem \ref{the10}, the Riccati equation (\ref{m5}) admits a solution $P>0.$
Since $B$ has a full column rank, then $B'PB+B'(A')^dPA^dB>0.$ Using the fact that $M^{-1}<N^{-1}$ when $M>N>0$ and $R>0,Q>0,$ we have
\begin{eqnarray}
P&>&A'PA-\gamma_iA'PB\Big[B'PB\nonumber\\
&&+B'(A')^dPA^dB\Big]^{-1}B'PA.\label{m26}
\end{eqnarray}
From $p\in(0,1),$ one has $\mu>0$ and $\sigma^2>0$ which yields that $\mu^2B'PB>0.$
Thus (\ref{m26}) further implies that
\begin{eqnarray}
P&>&A'PA-\bar{\gamma}_i(L_G)A'PB\Big[\mu^2B'PB\nonumber\\
&&+\sigma^2B'(A')^dPA^dB\Big]^{-1}B'PA,\label{m23}
\end{eqnarray}
where $\bar{\gamma}_i=\mu^2\frac{4\Big(\lambda_i(L_G)[\lambda_2(L_G)+\lambda_N(L_G)]-\lambda_i^2(L_G)\Big)}{[\lambda_N(L_G)+\lambda_2(L_G)]^2}.$
By letting the feedback gain matrix
\begin{eqnarray}
K&=&\frac{2\mu}{\lambda_2(L_G)+\lambda_N(L_G)}\Big[\mu^2B'PB\nonumber\\
&&+\sigma^2B'(A')^dPA^dB\Big]^{-1}B'PA,\label{m24}
\end{eqnarray}
the Riccati equation (\ref{m23}) is equivalently rewritten as
\begin{eqnarray}
P&>&[A-\lambda_i(L_G)\mu BK]'P[A-\lambda_i(L_G)\mu BK]\nonumber\\
&&+\sigma^2\lambda_i^2(L_G)K'B'(A')^dPA^dBK.\label{m6}
\end{eqnarray}
Combining with Lemma \ref{the1}, system (\ref{m4}) is mean-square stabilizable. This yields that the multi-agent system (\ref{m18}) is mean-square consensusable.
The proof is now completed. \hfill $\blacksquare$

\begin{remark}\label{rem1}
Noting that $\mu=1-p$ and $\sigma^2=p(1-p),$ the condition $\gamma_2>\gamma_c$ in Theorem \ref{the2} becomes
$(1-p)\Big[1-\Big(\frac{\lambda_N(L_G)-\lambda_2(L_G)}{\lambda_N(L_G)+\lambda_2(L_G)}\Big)^2\Big]>\gamma_c.$
\end{remark}

\begin{remark}

When time delay is reduced to $0,$ the sufficient condition
$\frac{\mu^2}{\mu^2+\sigma^2}\Big[1-\Big(\frac{\lambda_N(L_G)-\lambda_2(L_G)}{\lambda_N(L_G)+\lambda_2(L_G)}\Big)^2\Big]>\gamma_c$
is consistent with the result obtained in \cite{liangxu} for the consensusability of discrete-time linear multi-agent systems over analog
fading networks where $\mu$ and $\sigma^2$ are corresponding to the expectation and the covariance of identical channel fading .

\end{remark}

We next give a necessary condition for the mean-square consensusability of multi-agent system (\ref{m18}).
\begin{theorem}\label{the4}
Under Assumption \ref{a1}, \ref{a2} and $Rank(B)=1,$ the multi-agent system (\ref{m18}) is mean-square consensusable under protocol
(\ref{m19}) only if
\begin{eqnarray}
\Pi_i|\lambda_i^u(A)|^2<\Big(\frac{1+\lambda_2(L_G)/\lambda_N(L_G)}{1-\lambda_2(L_G)/\lambda_N(L_G)}\Big)^2,\label{m20}
\end{eqnarray}
where $\lambda_i^u(A)$ denotes the unstable eigenvalue of matrix $A.$
\end{theorem}
\emph{Proof.}
Using Theorem \ref{the3}, systems (\ref{m4}) are mean-square stable simultaneously for all $i=2,\ldots,N.$
By applying Lemma \ref{the1}, the following systems \begin{eqnarray}
\tilde{\delta}_i(k+1)&=&A\tilde{\delta}_i(k)-\lambda_i\mu BK\tilde{\delta}_i(k)-w_k \lambda_iA^dBK\tilde{\delta}_i(k)\nonumber
\end{eqnarray}
are mean-square stable for all $i=2,\ldots,N.$
Combining with the fact that $\lim_{k\rightarrow\infty}E\|\tilde{\delta}_i(k)\|^2=0$ implies that
$\lim_{k\rightarrow\infty}E\big[\tilde{\delta}_i(k)\big]=0,$ it yields that $A-\lambda_i\mu BK$ is Schur stable, i.e. all the eigenvalues
of $A-\lambda_i\mu BK$ are within the unit disk. The result then follows from \cite{you}. \hfill $\blacksquare$

\begin{remark}
Consider the case of $Rank(B)=1.$
When the communication is delay free and packets can be perfectly delivered, that is,
$d=0$ and $p=0,$ $\gamma_c=1-\frac{1}{\Pi_i|\lambda_i(A)|^2}$ which has been obtained in \cite{sinopoli1}.
From Theorem \ref{the2}, $\gamma_2>\gamma_c$ is reduced to (\ref{m20}). Together with Theorem \ref{the4},
(\ref{m20}) is necessary and sufficient for the consensusability of multi-agent systems (\ref{m18}) under protocol (\ref{m19}).
This is consistent with Theorem 3.1 in \cite{you} for the deterministic linear multi-agent systems under\begin{eqnarray}
u_i(k)=K\sum_{j\in N_i}\big[x_j(k)-x_i(k)\big].\nonumber
\end{eqnarray}
\end{remark}
%

We then study the scalar multi-agent systems. It shall be shown that $\gamma_2>\gamma_c$ in Theorem \ref{the2} is necessary and sufficient for the consensusability.
\begin{theorem}\label{the5}
Let $A=a\geq1, B=b>0 $ be constants, the multi-agent system (\ref{m18}) is mean-square consensusable by the control protocol (\ref{m19})
if and only if
\begin{eqnarray}
\frac{\mu^2}{(\mu^2+a^{2d}\sigma^2)}\Big[1-\frac{\Big(\lambda_N(L_G)-\lambda_2(L_G)\Big)^2}{\Big(\lambda_2(L_G)+\lambda_N(L_G)\Big)^2}\Big]> 1-\frac{1}{a^2}\label{m15}
\end{eqnarray}
\end{theorem}
\emph{Proof.} The equivalent condition (\ref{m6}) for the consensusability is reduced to
\begin{eqnarray}
&&\hspace{-6mm}a^2-2\lambda_i(L_G)\mu abk+\lambda_i^2(L_G)\mu^2b^2k^2+a^{2d}\sigma^2\lambda_i^2(L_G)b^2k^2\nonumber\\
&<&1,\label{m25}
\end{eqnarray}
that is,
\begin{eqnarray}
\lambda_i(L_G)^2(\mu^2+a^{2d}\sigma^2)b^2k^2-2\lambda_i(L_G)\mu abk+a^2-1<0.\label{m14}
\end{eqnarray}
``Necessity"
Since $b>0,$ one has from (\ref{m14}) that
\begin{eqnarray}
&&\frac{\mu a-\sqrt{(\mu a)^2-(\mu^2+a^{2d}\sigma^2)(a^2-1)}}{\lambda_i(L_G)(\mu^2+a^{2d}\sigma^2)b}\leq k\nonumber\\
&\leq&\frac{\mu a+\sqrt{(\mu a)^2-(\mu^2+a^{2d}\sigma^2)(a^2-1)}}{\lambda_i(L_G)(\mu^2+a^{2d}\sigma^2)b}\nonumber
\end{eqnarray}
Thus, we obtain that $\bigcap_{i=2}^N\Big(\frac{\mu a-\sqrt{(\mu a)^2-(\mu^2+a^{2d}\sigma^2)(a^2-1)}}{\lambda_i(L_G)(\mu^2+a^{2d}\sigma^2)b},\\
\frac{\mu a-\sqrt{(\mu a)^2-(\mu^2+a^{2d}\sigma^2)(a^2-1)}}{\lambda_i(L_G)(\mu^2+a^{2d}\sigma^2)b}\Big)\neq\emptyset.$ Using $\lambda_2(L_G)<\lambda_i(L_G)<\lambda_N(L_G),$ it
is further derived that
\begin{eqnarray}
&&\frac{\mu -\sqrt{(\mu )^2-(\mu^2+a^{2d}\sigma^2)(1-\frac{1}{a^2})}}{\lambda_2(L_G)}\nonumber\\
&\leq&\frac{\mu +\sqrt{(\mu )^2-(\mu^2+a^{2d}\sigma^2)(1-\frac{1}{a^2})}}{\lambda_N(L_G)}.\nonumber
\end{eqnarray}
By applying some algebraic transformations, we have
\begin{eqnarray}
\Big[\frac{\Big(\lambda_N(L_G)-\lambda_2(L_G)\Big)^2}{\Big(\lambda_2(L_G)+\lambda_N(L_G)\Big)^2}-1\Big]\mu^2
\leq -(\mu^2+a^{2d}\sigma^2)(1-\frac{1}{a^2}).\nonumber
\end{eqnarray}
Thus, (\ref{m15}) follows.

``Sufficiency" From (\ref{m15}), it yields that
\begin{eqnarray}
&&\hspace{-6mm}\frac{\mu^2}{(\mu^2+a^{2d}\sigma^2)}\frac{4\Big[\lambda_i(L_G)\Big(\lambda_N(L_G)+\lambda_2(L_G)\Big)-\lambda_i^2(L_G)\Big]}{(\lambda_2(L_G)+\lambda_N(L_G))^2}\nonumber\\
&>&1-\frac{1}{a^2}.\nonumber
\end{eqnarray}
Selecting the feedback gain in the form of (\ref{m24}) which gives that
$k=\frac{2\mu}{(\mu^2+\sigma^2 a^{2d})\Big(\lambda_2(L_G)+\lambda_N(L_G)\Big)}\frac{a}{b}.$
Then (\ref{m25}) follows. Thus, system (\ref{m18}) is mean-square consensusable.
The proof is now completed. \hfill $\blacksquare$

\begin{remark}
For system (\ref{m18}) with delay and $p=0$, the advantage of using the predictor-like protocol (\ref{m19}) is that the allowable delay for consensus can be arbitrarily large. However, when using the protocol without delay compensation, there exists a maximum delay margin within which consensus can be achieved \cite{huisa}. Take the case of $Rank(B)=1$ for example, by combining Theorem 3, Theorem 4 with Lemma 5.4 in \cite{sinopoli}, the equivalent condition for consensus of system (\ref{m18}) is $
\Pi_i|\lambda_i^u(A)|^2<\Big(\frac{1+\lambda_2(L_G)/\lambda_N(L_G)}{1-\lambda_2(L_G)/\lambda_N(L_G)}\Big)^2.$
This is exactly the necessary and sufficient condition to ensure the consensus for system (\ref{m18}) without delay obtained in \cite{you}. This indicates that system (\ref{m18}) is consensusable for any large delay under the basic assumption.
Furthermore, recalling Theorem 3 in \cite{xuzhangxie}, for scalar system with input delay, when $1+\frac{\lambda_2(L_G)}{\lambda_N(L_G)}\leq A<\frac{1+\lambda_2(L_G)/\lambda_N(L_G)}{1-\lambda_2(L_G)/\lambda_N(L_G)}$ or
    $-\frac{1+\lambda_2(L_G)/\lambda_N(L_G)}{1-\lambda_2(L_G)/\lambda_N(L_G)}\leq A\leq-1,$ no delay is allowed for consensusability via relative state feedback protocols. This illustrates the advantage of using predictor-like protocol (\ref{m19}) which can tolerate any large delay.
\end{remark}

As an important application, the result on consensusability is
extended to study formationability of the discrete-time multi-agent systems (\ref{m18}). In particular, given a formation vector
$H=\left[
     \begin{array}{ccc}
       H_1' & \cdots & H_N' \\
     \end{array}
   \right]'
$, the following control protocol is adopted
to study the formation problem of the discrete-time multi-agent
systems:
\begin{eqnarray}
u_i(k-d)&=&K\sum_{j\in N_i}\Big(\big[\hat{x}_j(k|k-d)-H_j\big]\nonumber\\
&&-\big[\hat{x}_i(k|k-d)-H_i\big]\Big),\label{m21}
\end{eqnarray}
where $H_i-H_j$ is the desired formuation vector between agent $i$
and agent $j.$ Noting that the common knowledge of the
directions of reference axes is required for all the agents, the protocol $u_i(k)=K\sum_{j\in N_i}\Big(\big[x_j(k)-H_j\big]-\big[x_i(k)-H_i\big]\Big)$ has been widely adopted in formation control \cite{you} and references therein, we now apply the predictor-like protocol (\ref{m19}) to the formationable problem.

\begin{definition}
The discrete-time multi-agent system (\ref{m18}) is said to be formationable under protocol
(\ref{m21}) if for any finite $x_i(0), u_i(-d),\cdots, u_i(-1),$ there exists a control gain $K$
in (\ref{m21}) such that $\lim_{k\rightarrow\infty}E\|\big[x_j(k)-H_j\big]-\big[x_i(k)-H_i\big]\|^2=0,\forall i, j=1,\cdots,N.$
\end{definition}
Based on Theorem \ref{the2}, sufficient and necessary conditions
on formationability of the discrete-time multi-agent systems is
stated as follows.

\begin{corollary}
Assume that Assumption \ref{a1} holds and $A(H_i-H_j)=(H_i-H_j), ~\forall i, j=1,\cdots,N.$ The following statements hold:
\begin{enumerate}
  \item If $\gamma_2>\gamma_c$ where $\gamma_c$ is given in Theorem \ref{the10},
then the multi-agent system (\ref{m18}) is mean-square formationable under protocol
(\ref{m21}).
  \item Let $Rank(B)=1,$ the multi-agent system (\ref{m18}) is mean-square consensusable under protocol
(\ref{m21}) only if (\ref{m20}) holds.
  \item Let $A=a\geq1,B=b>0,$ the multi-agent system (\ref{m1}) is mean-square formationable under protocol (\ref{m21})
if $\frac{\mu^2}{(\mu^2+a^{2d}\sigma^2)}\Big[1-\frac{\Big(\lambda_N(L_G)-\lambda_2(L_G)\Big)^2}{\Big(\lambda_2(L_G)+\lambda_N(L_G)\Big)^2}\Big]> 0.$
\end{enumerate}
\end{corollary}
\emph{Proof.} Denote $\delta_i(k)=\big[x_i(k)-H_i\big]-\big[\bar{X}(k)-\bar{H}\big]$ where $\bar{X}(k)=\frac{1}{N}\sum_{i=1}^Nx_i(k),$ $\bar{H}=\frac{1}{N}\sum_{i=1}^NH_i.$ Then mean-square formationability is equivalent to
that $\lim_{k\rightarrow\infty} E\|\delta_i(k)\|^2=0.$ By stacking $\delta_i$ into a column vector $\delta(k)=\left[
                                                                                                                \begin{array}{ccc}
                                                                                                                  \delta_1'(k) & \cdots & \delta_N'(k) \\
                                                                                                                \end{array}
                                                                                                              \right],
$ the following dynamical equation is in force:
\begin{eqnarray}
\delta(k+1)
&=&(I_N\otimes A)\delta(k)-\mu(L_G\otimes BK)\hat{\delta}(k|k-d)\nonumber\\
&&- w_k(L_G\otimes BK)\hat{\delta}(k|k-d)\nonumber\\
&&+\big[I_N\otimes (A-I_n)\big]\left[
                                                                        \begin{array}{c}
                                                                          H_1-\bar{H} \\
                                                                          \vdots \\
                                                                          H_N-\bar{H} \\
                                                                        \end{array}
                                                                      \right].\nonumber
\end{eqnarray}
Together with  $A(H_i-H_j)=(H_i-H_j),$ it follows that $(A-I_n)\big](H_1-\bar{H})=0.$ The above equation is thus reformulated as
\begin{eqnarray}
\delta(k+1)
&=&(I_N\otimes A)\delta(k)-\mu(L_G\otimes BK)\hat{\delta}(k|k-d)\nonumber\\
&&- w_k(L_G\otimes BK)\hat{\delta}(k|k-d).\nonumber
\end{eqnarray}
The remainder of the proof follows from Theorem \ref{the3}, \ref{the2}, \ref{the4} and \ref{the5}. The proof is now completed. \hfill $\blacksquare$

\section{Conclusions}

In this paper, we studied the consensusability of multi-agent systems with delay and packet dropout.
By proposing a kind of predictor-like protocol, sufficient and necessary conditions have been given for the mean-square
consensusability in terms of system matrices, time delay, communication graph and the packetdrop probability. It has been shown that the
derived results are exactly the necessary and sufficient condition obtained in \cite{you} for the delay and packet drop free.
Moreover, sufficient and necessary conditions have been obtained for the formationability of multi-agent systems.


\appendix

The following results can be obtained by similar discussions as in \cite{sinopoli1}. We give some brief proofs for the completion of the work.
\begin{lemma}\label{lem2}
Assume that $P\in \{S\in R^{n\times n}, S\geq0\},R>0, Q>0.$ Then the following statements hold.
\begin{enumerate}
  \item With $K_P=-\Big[R+B'PB+B'(A')^dPA^dB\Big]^{-1}B'PA,~g_{\gamma}(P)=\Phi(K_P,P).$
  \item $g_{\gamma}(P)=\min_{K}\Phi(K,P)\leq \Phi(K,P).$
  \item If $P_1\leq P_2,$ then $g_{\gamma}(P_1)\leq g_{\gamma}(P_2).$
  \item If $\gamma_1\leq \gamma_2,$ then $g_{\gamma_1}(P)\geq g_{\gamma_2}(P).$
  \item If $\alpha\in [0,1],$ then $g_{\gamma}\Big(\alpha P_1+(1-\alpha)P_2\Big)\geq \alpha g_{\gamma}(P_1)+(1-\alpha)g_{\gamma}(P_2).$
  \item $g_{\gamma}(P)\geq (1-\gamma)A'PA+Q.$
  \item Provided that the equation $(1-\gamma)A'XA+Q=X$ has a solution $X>0.$ If $\bar{P}\geq g_\gamma(\bar{P}),$ then $\bar{P}>0$
\end{enumerate}
\end{lemma}
\emph{Proof.}
\begin{enumerate}
  \item Using the definition of $K_P,$ we have
\begin{eqnarray}
\Phi(K_P,P)&=&A'PA+Q-\gamma A'PB \Big[R+B'PB\nonumber\\
&&+B'(A')^dPA^dB\Big]^{-1}B'PA\nonumber\\
&=&g_{\gamma}(P).\nonumber
\end{eqnarray}
  \item By using the definitions of $\Phi(K,P)$ and $\Psi(K,P),$ it holds that
$\min_{K}\Phi(K,P)=\min_{K}\Psi(K,P).$
Combining with the fact that $P\geq0, R>0,$ the minimum of $K$ can be found by using $\frac{\partial\Psi(K,P) }{\partial K}=0,$ that is
$0=B'P(A+BK)+B'(A^d)'RA^dBK+RK.$
This implies that
$K=-\Big[R+B'PB+B'(A')^dPA^dB\Big]^{-1}B'PA.$
Together with from fact 1), the result follows.
  \item If $P_1\leq P_2,$ we have by using the above two facts
\begin{eqnarray}
g_{\gamma}(P_1)&=&\Phi(K_{P_1},P_1)\leq\Phi(K_{P_2},P_1)\nonumber\\
 &\leq&\Phi(K_{P_2},P_2)= g_{\gamma}(P_2).\nonumber
\end{eqnarray}
  \item Noting that $A'PB\Big[R+B'PB+B'(A')^dPA^dB\Big]^{-1}B'PA\geq0,$ the fact follows directly.
  \item Let $Z=\alpha P_1+(1-\alpha) P_2,$ then
\begin{eqnarray}
g_{\gamma}(Z)
&=&(1-\gamma)(A'ZA+Q)+\gamma \Psi(K_Z,Z).\nonumber
\end{eqnarray}
Further rewriting $\Psi(K_Z,Z)$ yields that
\begin{eqnarray}
\Psi(K_Z,Z)
&=&\alpha \Psi(K_Z,P_1)+(1-\alpha)\Psi(K_Z,P_2)\nonumber\\
&\geq& \alpha \Psi(K_{P_1},P_1)+(1-\alpha)\Psi(K_{P_2},P_2)\nonumber.
\end{eqnarray}
Thus
\begin{eqnarray}
g_{\gamma}(Z)
&\geq&(1-\gamma)(A'ZA+Q)+\gamma \alpha \Psi(K_{P_1},P_1)\nonumber\\
&&+\gamma (1-\alpha)\Psi(K_{P_2},P_2)\nonumber\\
&=&\alpha g_{\gamma}(P_1)+(1-\alpha)g_{\gamma}(P_2).\nonumber
\end{eqnarray}
  \item By using the facts that $F_1'PF_1\geq0, F_2'PF_2\geq0, K'RK\geq0,$ the result is straightforward.
  \item Using the above fact, it follows that $\bar{P}\geq g_{\gamma}(\bar{P})\geq (1-\gamma)A'\bar{P}A+Q.$ Combining with $(1-\gamma)A'XA+Q=X,$
 there holds that $\bar{P}-X\geq (1-\gamma)A'(\bar{P}-X)A,$ which gives $\bar{P}-X\geq0.$ Since $X>0,$ it is thus obtained that $\bar{P}>0.$
\end{enumerate}

\begin{theorem}\label{the6}
Suppose there exists a matrix $\tilde{K}$ and a positive-definite matrix $\tilde{P}$ such that $\tilde{P}>\Phi(\tilde{K},\tilde{P}).$ Then
\begin{enumerate}
  \item for any initial condition $P_0$, the MARE converges,
and the limit is independent of the initial condition
$\lim_{t\rightarrow\infty}P_t=\lim_{t\rightarrow\infty}g_\gamma^t(P_0)=\bar{P}.$
  \item $\bar{P}$ is the unique positive-semidefinite fixed point of the
MARE.
\end{enumerate}
\end{theorem}
\emph{Proof.}
\begin{enumerate}
  \item We first let the initial condition be $Q_0=0.$ Let $Q_k=g_{\lambda}^k(0).$ Since $0=Q_0\leq Q_1=Q.$ From $3)$ of Lemma \ref{lem2},
it follows that $Q_1=g_{\gamma}(Q_0)\leq g_{\gamma}(Q_1)=Q_2.$ By induction, it is obtained that $Q_t\leq Q_{t+1}$ for $t\geq0.$ We show the sequence
has an upper bound.
Define the linear operator
$\mathcal{L}(Y)=(1-\gamma)A'YA+\gamma(F_1'YF_1+F_2'YF_2).$
Noting that $\tilde{P}>\Phi(\tilde{K},\tilde{P})=\mathcal{L}(\tilde{P})+Q+\gamma K'RK\geq \mathcal{L}(\tilde{P}).$
On the other hand, we have $Q_{t+1}=g_{\gamma}(Q_t)\leq \Phi(K_{\tilde{P}},Q_t)=\mathcal{L}(\tilde{P})+Q+\gamma K_{\tilde{P}}'RK_{\tilde{P}}.$ In view of
$Q+\gamma K_{\tilde{P}}'RK_{\tilde{P}}\geq0$ and using Lemma 3 in \cite{sinopoli1}, we conclude that there exists $M_{Q_0}$ such that $Q_t\leq M_{Q_0}$ for $t\geq 0.$ Accordingly, the sequence converges, i.e. $\lim_{t\rightarrow\infty}Q_t=\bar{P}$ and $\bar{P}=g_{\gamma}(\bar{P}).$

We next consider the case that the initial condition is selected as $R_0\geq \bar{P}.$ First, define
$\bar{K}=-\Big[R+B'\bar{P}B+B'(A')^d\bar{P}A^dB\Big]^{-1}B'\bar{P}A,\bar{F}_1=A+B\bar{K},\bar{F}_2=A^dB\bar{K}$ and
$\hat{\mathcal{L}}(Y)=(1-\gamma)A'YA+\gamma(\bar{F}_1'Y\bar{F}_1+\bar{F}_2'Y\bar{F}_2).$
It is noted that
$\bar{P}=g_{\gamma}(\bar{P})=\hat{\mathcal{L}}(Y)+Q+\bar{K}'R\bar{K}>\hat{\mathcal{L}}(Y)$
where $Q>0$ has been used in the derivation of last inequality.
Using again Lemma 3 in \cite{sinopoli1}, we have that $\lim_{t\rightarrow\infty}\hat{\mathcal{L}}^t(Y)=0$ for all $Y\geq0.$
Since $R_0\geq \bar{P},$ then $R_1=g_{\gamma}(R_0)\geq g_{\gamma}(\bar{P})=\bar{P}.$ By induction, it follows that $R_t\geq \bar{P}$ for $t\geq0.$
Noting that
\begin{eqnarray}
0&\leq &R_{t+1}-\bar{P}=g_{\gamma}(R_t)-g_{\gamma}(\bar{P})\nonumber\\
&=&\Phi(K_{R_t},R_t)-\Phi(K_{\bar{P}},\bar{P})\nonumber\\
&\leq& \Phi(K_{\bar{P}},R_t)-\Phi(K_{\bar{P}},\bar{P})\nonumber\\
&=&(1-\gamma)A'(R_{t}-\bar{P})A+\gamma \bar{F}_1'(R_{t}-\bar{P})\bar{F}_1\nonumber\\
&&+\gamma\bar{F}_2'(R_{t}-\bar{P})\bar{F}_2=\hat{\mathcal{L}}(R_{t}-\bar{P})\rightarrow0, t\rightarrow\infty,\nonumber
\end{eqnarray}
which gives that
$\lim_{t\rightarrow\infty}R_{t+1}=\bar{P}.$

We now prove that the Riccati iteration converges to $\bar{P}$ for all initial values $P_0\geq0.$ Let $Q_0=0$ and $R_0=P_0+\bar{P},$ it is obvious
that $Q_0\leq P_0\leq R_0.$ Consider the Riccati iterations initialized at $Q_0, P_0 $ and $R_0.$ It then follows that
$Q_t\leq P_t\leq R_t, \forall t\geq0.$
Based on the above discussions, it has already been obtained that $\lim_{t\rightarrow\infty}Q_t=\lim_{t\rightarrow\infty}R_t=\bar{P}.$ This implies that
$\lim_{t\rightarrow\infty}P_t=\bar{P}.$

\item It is now claimed that the solution is unique. Otherwise, let $\hat{P}$ be another solution, i.e., $\hat{P}=g_{\gamma}(\hat{P})$ and let the initial value
be $\hat{P}.$ Thus we have a constant sequence with $\hat{P}.$ Using the above prove, we have that the constant sequence also converges to $\bar{P}.$
Thus $\hat{P}=\bar{P}.$
The proof is now completed.

\end{enumerate}

%
%
%
%

\begin{theorem}\label{the7}
If $(A,B,0,A^dB)$ is mean-square stabilizable and $A$ is unstable. Then there exists a $\gamma_c\in [0,1)$ such that
\begin{eqnarray}
\lim_{t\rightarrow\infty}P_t=+\infty, \mbox{for}~0\leq\gamma\leq \lambda_c~\mbox{and}~ \exists P_0\geq 0\nonumber\\
P_t\leq M_{P_0}~\forall t, \mbox{for}~\lambda_c<\gamma\leq 1~\mbox{and}~ \forall P_0\geq 0\nonumber
\end{eqnarray}
where $M_{P_0}>0$ depends on the initial condition $P_0\geq0.$
\end{theorem}
\emph{Proof.} If $\lambda=1,$ the Riccati difference equation becomes the delay-dependent Riccati equation in \cite{zhang1} and \cite{zhang3}
which has been shown to converge to a unique positive definite solution under the mean-square stabilizability of $(A,B,0,A^dB)$ for the zero initial value.
Based on similar discussions in Theorem \ref{the6}, the Riccati iteration converges to a fixed point for any initial values $P_0\geq0.$ Hence, $P_t$ is always bounded
for any initial values $P_0\geq0.$
If $\lambda=0,$ the equation is reduced to $P_{t+1}=A'P_tA+Q.$ If $A$ is unstable, there always exists one initial value$P_0\geq0$ such that
$P_t$ is unbounded. Accordingly, the critical value $\lambda_c\in [0,1)$ exists. We now prove there exists a single critical value. In fact, for any $\lambda>\lambda_c,$
it is obtained that $P_{t+1}=g_{\lambda}(P_t)\leq g_{\lambda_c}(P_t)$ which is bounded. This completes the proof.

\begin{theorem}\label{the8}
If $(A,B,0,A^dB)$ is mean-square stabilizable and $A$ is unstable. Then the critical value satisfies
$\underline{\gamma}\leq \gamma_c\leq \overline{\gamma}$ where
\begin{eqnarray}
\underline{\gamma}&=&arginf_{\gamma}\{\exists S|(1-\lambda)A'SA+Q=S, S\geq0\}\nonumber\\
\overline{\gamma}&=&arginf_{\gamma}\{\exists (K,P)|P>\Phi(K,P)\}\nonumber
\end{eqnarray}
\end{theorem}
\emph{Proof.} Consider $S_{t+1}=(1-\gamma)A'S_tA+Q$ with $S_0=0,$ it is obtained that $\lim_{t\rightarrow\infty}S_t=\infty$ for $\lambda>\underline{\gamma}$ in the proof of Theorem 3 in \cite{sinopoli1}. Noting that the initial value $P_0\geq0,$ i.e. $P_0\geq S_0.$ Assume that $P_t\geq S_t.$ From 6) of Lemma \ref{lem2}, it holds that
$P_{t+1}\geq (1-\gamma)A'P_tA+Q\geq(1-\gamma)A'S_tA+Q=S_{t+1}.$
By induction, we have that $P_t\geq S_t,\forall t\geq 0, \forall P_0\geq0.$ This implies that
$\lim_{t\rightarrow\infty}P_t\geq \lim_{t\rightarrow\infty}S_t=\infty.$ That is, $P_t$ is unbounded for any $\gamma<\underline{\gamma}$ and any initial values $P_0\geq0.$ Therefore, $\gamma_c\geq \underline{\gamma}.$
On the other hand, when $\gamma>\overline{\gamma},$ there exists $X$ such that $X>\Phi(K,X)\geq g_{\gamma}(X).$ Using 7) of Lemma \ref{lem2}, it yields that
$X>0.$ Using Lemma 3 of \cite{sinopoli1}, $P_t$ is bounded. That is, $\gamma_c\leq\underline{\gamma}.$

\begin{theorem}\label{the9}
If $(A,B,0,A^dB)$ is mean-square stabilizable, then the following statements are equivalent.
\begin{enumerate}
  \item $\exists X$ such that $X>g_{\gamma}(X).$
  \item $\exists K,X>0$ such that $X>\Phi(K,X).$
  \item $\exists Z$ and $0\leq Y\leq I$ such that
\begin{eqnarray}
\Gamma_\gamma(Y,Z)&=&\left[
                      \begin{array}{cccc}
                        Y & \sqrt{\gamma}(AY+BZ)' &  &  \\
                        \sqrt{\gamma}(AY+BZ) & Y &  &  \\
                        \sqrt{\gamma}A^dBZ& 0 &  &  \\
                        \sqrt{1-\gamma}AY & 0 &  &  \\
                      \end{array}
                    \right.\nonumber\\
&&\hspace{-6mm}\left.
                      \begin{array}{cccc}
                         &  & \sqrt{\gamma}(A^dBZ)' & \sqrt{1-\gamma}(AY)' \\
                         &  & 0 & 0 \\
                        &  & Y & 0 \\
                        &  & 0 & Y \\
                      \end{array}
                    \right]>0.\nonumber
\end{eqnarray}
\end{enumerate}
\end{theorem}
\emph{Proof.} Using facts 1) and 2) in Lemma \ref{lem2}, the equivalence between 1) and 2) follows.
We now establish the equivalence between 2) and 3). Let $F=A+BK,$ then $X>\Phi(K,X)$ is in fact
$X>(1-\gamma)A'XA+\gamma F'XF+\gamma K'B'(A')^dXA^dBK+\gamma K'RK+Q.$
By using Schur complement, the inequality is equivalent to
\begin{eqnarray}
\left[
  \begin{array}{cc}
    X-(1-\gamma)A'XA+\lambda K'B'(A')^dXA^dBK & \sqrt{\gamma} F' \\
     \sqrt{\gamma} F& X^{-1} \\
  \end{array}
\right]>0.\nonumber
\end{eqnarray}
By taking similar procedures to Theorem 5 in \cite{sinopoli1}, the result can be obtained. So we omit the details.


\end{document}